\documentstyle{article}

\newtheorem{theorem}{Theorem}

\begin{document}

\title{Fundamental group for some cuspidal curves \footnote{1991 
Mathematics Subject Classification 14H20, 14H30, 14E20}}

\author{ Jos\'e Ignacio Cogolludo \footnote{This work was 
partially supported by CAICYT PB94-0291}}

\maketitle

\section*{Introduction}

In \cite{hirano}, Hirano gives a method to construct families of curves
with a large number of singularities. The idea is to consider an abelian 
covering of ${\bf P}^2$ ramified along three lines in general position 
and to take the pull-back of a curve $C$ intersecting the lines 
non-generically. Similar constructions are used by Shimada in 
\cite{shim} and Oka in \cite{oka:two}. We apply this method for the case 
where $C$ is a conic, 
constructing a family of curves with the following asymptotic behavior 
(see \cite{zaid})
$$\lim_{n \rightarrow \infty} \Big({\sum_{p \in Sing \tilde{C}_n} \mu (p)}\Big)
\Big/ \Big({\deg(\tilde{C}_n)}\Big)^2 = \frac{3}{4}.$$

The goal of this paper is to calculate the fundamental group for the curves 
in this family as well as their Alexander polynomial. We have the following

\begin{theorem}
Let $f(x,y,z)=x^2+y^2+z^2+2(xz-xy+yz)$ and 
$\tilde{C}_n$ be the projective plane curve defined by 
$\{ {[}x:y:z{]} \in {\bf P}^2 \mid f(x^n,y^n,z^n)=0 \}$ then 
$$\array{lll}
 & \varepsilon_{i+2,j}\ =\  \varepsilon_{i+1,j}^{-1} \cdot \varepsilon_{i,j} 
\cdot \varepsilon_{i+1,j}, & \\
\pi_1({\bf P}^2 \setminus \tilde{C}_n)\ = \ < 
\varepsilon_{i,j}\ ;\ & 
\varepsilon_{i,j+2}\ =\ \varepsilon_{i,j+1}^{-1} \cdot \varepsilon_{i,j} \cdot
\varepsilon_{i,j+1},
\ \ \ \ \ \ 
i,j \in {\bf Z}/n{\bf Z} & >. \\ 
 & \varepsilon_{0,0} \cdot \varepsilon_{0,1} \cdot \varepsilon_{1,1}
\cdot ... \cdot \varepsilon_{n-1,n-1} \cdot \varepsilon_{n-1,0}=1 &
\endarray$$
\end{theorem}

All curves $\tilde{C}_{2k}$ are reducible.
The first irreducible member of this family, that is $\tilde{C}_3$, is the
sextic with $9$ ordinary cusps considered by Zariski in \cite{zar:genusp}.
He gave a presentation of its fundamental group (also see \cite{kan}). We 
show explicitly an isomorphism between ours and Zariski's presentation. 
The global Alexander polynomial of this sextic, $(t^2-t+1)^3$ (see 
\cite{libg:homot}), is different from the local Alexander polynomial of 
its singularities, $(t^2-t+1)$. This was the first and unique 
example in the literature. Our family provides infinitely 
many examples of this behavior for curves containing any type of 
locally irreducible $A_{2k}$-singularity. We have the following

\begin{theorem}
Let $\tilde{C}_n$ be a curve defined as above, 
then the global Alexander polynomial of $\tilde{C}_{2k+1}$ is
$$\Delta_{\tilde{C}_{2k+1}}(t)=(t^{2k}-t^{2k-1}+...+t^2-t+1)^3.$$
\end{theorem}

For $n$ odd, One can also relate the groups $\pi_1({\bf P}^2 \setminus 
\tilde{C}_n)$ with the free products $({\bf Z}/2{\bf Z}) \ast 
({\bf Z}/n{\bf Z})$ considered by Oka in \cite{oka:nonab}. The latter
correspond to fundamental groups of the curves defined by 
$\{ {[}x:y:z{]} \in {\bf P}^2 \mid (y^n-z^n)^2+(x^2-y^2)^n=0\}$. 
Our groups are non-central extensions of these.
In fact, we can obtain the Oka curves as a perturbation of ours. 

For the case $n=3$, this perturbation provides sextics with six ordinary 
cusps on a conic. The special position of their singularities results 
in the following, these sextics and those with six cusps in general 
position belong to different connected components, as Zariski showed 
in \cite{zar2}. Our construction provides a degeneration from the family 
of sextics with six ordinary cusps on a conic to the family of sextics 
with nine cusps. The latter family is connected, since any curve in it 
is dual to a non-singular cubic.

A more detailed summary of this article is the following: In the 
first section we apply the Hirano method for a conic $C$ tangent to 
three lines, obtaining a family of curves with $A_{n}$-singularities where 
$n$ depends on the degree of the covering. In the second section we calculate 
the fundamental group of the curves of this family and in the third section 
we find their Alexander polynomials. In section four we show how this family 
can be obtained also as a degeneration of the Oka curves described above and 
how the change of the fundamental groups can be detected by deforming the
original curve $C$.

\subsection*{Acknowledgements}
This work would not have been possible without the advice of Prof. Anatoly 
Libgober. Thanks for your interest and support.

\section{Construction of the family of curves $\tilde{C}_n$}
\label{sec-construction}
Let's consider three lines $L_0$, $L_1$ and $L_2$ in general position and
a non-singular conic $C\equiv\{f=0\}$ in ${\bf P}^2$ tangent to the lines.
By changing coordinates on the plane we can assume $L_i$ to 
be the line $\{ [x_0:x_1:x_2] \mid x_i=0 \}$, we define a map on 
${\bf P}^2$ as follows
$$\array{cccc}
[ x_0:x_1:x_2 ] & \in & {\bf P}^2 & \\
\downarrow  & & \downarrow & p_n\\
{[}x_0^n:x_1^n:x_2^n{]} & \in & {\bf P}^2 & \\
\endarray$$
for any natural number $n$.

This is an abelian regular covering outside the lines $L_0$, $L_1$ and $L_2$.
Its Galois group is $({\bf Z}/n{\bf Z}) \oplus ({\bf Z}/n{\bf Z})$.
On the lines it has degree $n$ except on the double points, where the map 
is one-to-one. Let's denote by $\tilde{L}_i$ the pre-image of the line $L_i$.
Now consider the pull-back of $C$, which is an algebraic curve defined by 
the equation $\tilde{C}_n \equiv \{f(x_0^n:x_1^n:x_2^n)=0\}$. Its degree is 
$2n$ and it has the following properties:

\begin{enumerate}
\parindent=1cm
\item[(1)]
{\em It has $3n$ singular points of local type $A_{n-1}$ and no other 
singularities.}

{\bf Proof.} \quad
Since the covering is unramified outside the three 
lines $\tilde{L}_0, \tilde{L}_1$ and 
$\tilde{L}_2$, the curve $\tilde{C}_n$ is non-singular outside this locus. 
Therefore it can only be singular on the intersection with the pre-image
of the branching locus. This corresponds to the pre-image of the three
points of intersection of $C$ and $L_0 \cup L_1 \cup L_2$, say $p_0$, $p_1$
and $p_2$. Indeed this pre-image consists of $3n$ points. The equation of 
$C$ around any of the point $p_i$ can be written locally as $y^2=x$, where
$\{x=0\}$ is the local equation of the branching line $L_i$. Since 
the ramification order is $n$, the local equation for the pre-images 
is $y^2=x^n$, that is, singular points of type $A_{n-1}$. \quad $\Box$

\item[(2)]
{\em If $n$ is odd $\tilde{C}_n$ is irreducible.}

{\bf Proof} \quad
Since all singularities of $\tilde{C}_n$ are locally 
irreducible, the curve must be irreducible. \quad $\Box$

\item[(3)]
{\em If $n$ is even $\tilde{C}_n$ has exactly four 
non-singular irreducible components.}

{\bf Proof} \quad 
If $n=2$ it turns out that $\tilde{C}_n$ has degree $4$ 
and $12$ singular points of type $A_1$, that is, ordinary double points. 
By the genus formula $\tilde{C}_n$ cannot be irreducible, and by Bezout 
theorem, $\tilde{C}_n$ cannot be the union of two non-singular conics. 
Similar considerations eliminate the cases of $\tilde{C}_n$ being a line 
and an irreducible cubic or a conic and two lines. Therefore, it must be 
an arrangement of lines intersecting transversely.

If $n=2k$, the covering factors through $p_n=\tilde{p}_k \circ p_2.$
Therefore $\tilde{C}_n$ will be the pre-image by the abelian covering
$\tilde{p}_k$ of four lines intersecting, transversely, the ramification
locus. Since transverse intersections with the ramification locus do not
generate singularities on their pre-images, the pull-back of each line 
is a non-singular curve and therefore irreducible. \quad $\Box$
\end{enumerate}

\section{Fundamental group of ${\bf P}^2 \setminus \tilde{C}_n$}
The first step will be to calculate a presentation for 
$$G:=\pi_1 ({\bf P}^2 \setminus (C \cup L_0 \cup L_1 \cup L_2 )).$$

Since we have a real picture of this arrangement where all singular 
points appear, we can proceed rather easily by using the Zariski-Van
Kampen method. To start with, we consider the line $L$ of Fig. $1$.

\vspace*{20 pt}
\begin{picture}(300,120)(0,0)
\centering
\put(150,60){\circle{40}}
\put(78,17){\line(1,1){88}}
\put(73,12){\makebox(0,0){$L_1$}}
\put(150,106){\line(1,-2){45}}
\put(200,10){\makebox(0,0){$L_0$}}
\put(72,40){\line(1,0){133}}
\put(65,40){\makebox(0,0){$L_2$}}
\thinlines
\put(101,10){\line(0,1){100}}
\put(101,115){\makebox(0,0){$i_1$}}
\put(129.5,10){\line(0,1){98}}
\put(128,113){\makebox(0,0){$i_2$}}
\put(136,10){\line(0,1){100}}
\put(136,115){\makebox(0,0){$i_3$}}
\put(150,10){\line(0,1){98}}
\put(148,113){\makebox(0,0){$i_4$}}
\put(155.6,10){\line(0,1){100}}
\put(156,115){\makebox(0,0){$i_5$}}
\put(168.5,10){\line(0,1){98}}
\put(166,113){\makebox(0,0){$i_6$}}
\put(170.5,10){\line(0,1){100}}
\put(174,115){\makebox(0,0){$i_7$}}
\put(183,10){\line(0,1){98}}
\put(186,113){\makebox(0,0){$i_8$}}
\thicklines
\put(117,5){\line(0,1){115}}
\put(117,0){\makebox(0,0){$L$}} 
\put(145,-10){\makebox(0,0){\small \rm{Fig. $1$}}}
\end{picture}
\vspace*{10pt}

\vspace*{20pt}
\begin{picture}(300,100)(0,0)
\centering
\thicklines
\put(0,0){\line(1,0){200}}
\put(100,100){\line(1,0){200}}
\put(0,0){\line(1,1){100}}
\put(200,0){\line(1,1){100}}
\put(0,50){\makebox(0,0){down $\leftarrow$}}
\put(275,50){\makebox(0,0){$\rightarrow$ up}}
\put(75,50){\circle*{3}}
\put(90,50){\makebox(0,0){$\ell_2$}}
\put(100,25){\circle*{3}}
\put(115,25){\makebox(0,0){$e_2$}}
\put(150,75){\circle*{3}}
\put(135,75){\makebox(0,0){$e_1$}}
\put(175,50){\circle*{3}}
\put(190,50){\makebox(0,0){$\ell_1$}}
\put(225,50){\circle*{3}}
\put(240,50){\makebox(0,0){$\ell_0$}}
\put(250,20){\makebox(0,0){$L$}}
\thinlines
\put(75,50){\oval(16,12)}                      
\put(100,25){\oval(16,12)}
\put(150,75){\oval(16,12)}
\put(175,50){\oval(16,12)}
\put(225,50){\oval(16,12)}
\put(77,44){\vector(1,0){0}}
\put(102,19){\vector(1,0){0}}
\put(152,69){\vector(1,0){0}}
\put(177,44){\vector(1,0){0}}
\put(227,44){\vector(1,0){0}}
\put(80,55){\line(1,1){40}}
\put(180,55){\line(1,1){20}} 
\put(230,55){\line(1,1){10}} 
\put(155,80){\line(1,1){10}}
\put(105,30){\line(3,2){78}}
\put(120,95){\line(1,0){175}}
\put(200,75){\line(1,0){75}}
\put(240,65){\line(1,0){25}}
\put(165,90){\line(1,0){125}}
\put(183,82){\line(1,0){99}}
\put(145,-20){\makebox(0,0){\small {\rm Fig. $2$}}}
\end{picture}

\vspace*{40 pt}

The paths $e_1, e_2, \ell_0, \ell_1$ and $\ell_2$ of Fig. $2$
correspond to a system of generators for $G$. A system of relators 
is given by the action of the monodromy around the exceptional lines
$i_1,...,i_8$. There is a canonical way of considering this action 
which comes from the fact that the arrangement is real 
\cite{mois:braid}.
Proceeding in this way, the relations look as follows:
$$
\displaylines{
e_1=e_2=e \quad ({\rm around \ }i_2); \cr
\rm{[} e\ell_1e^{-1},\ell_2 \rm{]} =1 \quad ({\rm around \ }i_1); \cr
(\ell_1e)^2=(e\ell_1)^2 \quad ({\rm around \ }i_3); \cr
(\ell_2e)^2=(e\ell_2)^2  \quad ({\rm around \ }i_4); \cr
\rm{[} e\ell_1e^{-1},\ell_0 \rm{]} =1 \quad ({\rm around \ }i_5); \cr
(\ell_0 \cdot \ell_1^{-1}e\ell_1)^2=
(\ell_1^{-1} e \ell_1 \cdot
\ell_0)^2 \quad ({\rm around \ }i_6); \cr
\ell_0^{-1}\ell_1^{-1}e\ell_1\ell_0=\ell_2e\ell_2^{-1} 
\quad ({\rm around \ }i_7); \cr
\rm{[} \ell_2, \ell_1^{-1}e\ell_1\ell_0 
\ell_1^{-1}e^{-1}\ell_1 \rm{]}=1 \quad ({\rm around \ }i_8); \cr
\ell_2e^2\ell_1\ell_0=1 \quad ({\rm from \ the \ projective \ situation.}) 
\cr}
$$

This presentation can be simplified to get the following
$$\array{lll}
& [e\ell_1e^{-1},\ell_2]=1, & \\
G \ = \ <\ e,\ell_1,\ell_2\ ;\ & (\ell_1e)^2=(e\ell_1)^2, \ (\ell_2e)^2=
(e\ell_2)^2, & \ >. \\
& (\ell_1\ell_2e)^2=(e\ell_1\ell_2)^2 &
\endarray$$

Since $p_{n \mid {\bf P}^2 \setminus (\tilde{L}_0 \cup \tilde{L}_1 
\cup \tilde{L}_2)}$ is a topological covering, we have the 
following short exact sequence
$$ 1 \rightarrow K \rightarrow G \rightarrow S = {\bf Z}/n{\bf Z}
\oplus {\bf Z}/n{\bf Z} \rightarrow 1 ,$$
where $S$ is isomorphic to the commutative subgroup of $S_{n^2}$ 
generated by the images of $\ell_1$ and $\ell_2$, and $K$ is the 
fundamental group of ${\bf P}^2 \setminus (\tilde{C} \cup \tilde{L}_0 
\cup \tilde{L}_1 \cup \tilde{L}_2).$

Our group $\pi_1({\bf P}^2 \setminus \tilde{C}_n)$ is
the quotient of $K$ by the normal subgroup generated by the
meridians of the three lines, that is, $\ell_1^n$, $\ell_2^n$ and 
$\ell_0^n=(\ell_2e^2\ell_1)^{-n}$. In order to make things easier, 
and since we are not interested in $K$, but in a quotient, we will 
use these new relations in the coming calculations. 

To use Reidemeister-Schreier method \cite{kar:comb},
we need a set of representatives for cosets of $G$ with respect 
to $K$, say $\{ \ell_1^i \ell_2^j \}_{i,j = 0,...,n-1}$, which gives 
rise to a system of generators for $K$ 
$$\array{lll}
e_{i,j} & := & (\ell_1^i \ell_2^j) \cdot e \cdot (\ell_1^i \ell_2^j)^{-1} \\
\ell_{1,i,j} & := & (\ell_1^i \ell_2^j) \cdot \ell_1 
\cdot (\ell_1^{i+1} \ell_2^j)^{-1} \\
\ell_{2,i} & := & \ell_1^i \cdot \ell_2^n \cdot \ell_1^{-i},
\endarray$$
where $i,j \in {\bf Z}/n{\bf Z}.$

Observe that the meridians of the lines $\tilde{L}_1$ and 
$\tilde{L}_2$ become $\ell_{1,n-1,0}$ and $\ell_{2,i}$. Therefore we can 
eliminate these generators whenever they appear in calculations.

To simplify the notation, let's denote
$$\array{ll}
e_{i,j}^{r,s}:=e_{r,s}^{-1} \cdot e_{i,j} \cdot e_{r,s}, & \ \  {\rm and} \\
\varepsilon_{i,j}:= e_{i,j}^{i,0}. & 
\endarray$$

In order to obtain a system of relators for $K$ we use the rewriting method 
for the relators of $G$. So, for $[e\ell_1e^{-1},\ell_2]=1$ we get
$$e_{i,j} \cdot \ell_{1,i,j} \cdot e_{i+1,j}^{-1} \cdot e_{i+1,j+1} \cdot
\ell_{1,i,j+1} \cdot e_{i,j+1}^{-1}=1,$$
which allows us to eliminate all generators of the form $\ell_{1,i,j}$ by

\begin{equation}
\label{eles}
\ell_{1,i,j}=e_{i,j}^{-1} \cdot e_{i,0}
\cdot e_{i+1,0}^{-1} \cdot e_{i+1,j}.
\end{equation}
Using (\ref{eles}) we can rewrite $(\ell_1e)^2=(e\ell_1)^2$ as
$$\varepsilon_{i+2,j}\ =\  \varepsilon_{i,j}^{i+1,j},$$
and $(\ell_2e)^2\ =\ (e\ell_2)^2$ as 
$$e_{i,j+2}\ =\ e_{i,j}^{i,j+1},$$ or equivalently
$$\varepsilon_{i,j+2}\ =\ \varepsilon_{i,j}^{i,j+1}.$$
The relations coming from $(\ell_1\ell_2e)^2=(e\ell_1\ell_2)^2$ are all 
consequences of the previous ones.

It only remains to add the last relation for the meridian 
$\ell_0^n$, which can be written as

$$\varepsilon_{0,0} \cdot \varepsilon_{0,1} \cdot \varepsilon_{1,1} 
\cdot ... \cdot \varepsilon_{n-1,n-1} \cdot \varepsilon_{n-1,0}=1.$$

Therefore, we have obtained the desired presentation

\begin{equation}
\label{pres}
\array{lll}
 & \varepsilon_{i+2,j}\ =\  \varepsilon_{i,j}^{i+1,j}, & \\
\pi_1({\bf P}^2 \setminus \tilde{C}_n)\ = \ < 
\varepsilon_{i,j}\ ;\ & 
\varepsilon_{i,j+2}\ =\ \varepsilon_{i,j}^{i,j+1},
\ \ \ \ \ \ 
i,j \in {\bf Z}/n{\bf Z} & >. \\ 
 & \varepsilon_{0,0} \cdot \varepsilon_{0,1} \cdot \varepsilon_{1,1}
\cdot ... \cdot \varepsilon_{n-1,n-1} \cdot \varepsilon_{n-1,0}=1 &
\endarray
\end{equation}

Note that it is fairly clear, by (\ref{pres}), how the 
parity of $n$ determines that the abelianization of $\pi_1({\bf P}^2 
\setminus \tilde{C}_n)$ is ${\bf Z}/2n{\bf Z}$ or ${\bf Z}^3 \oplus
{\bf Z}/\frac{n}{2}{\bf Z}$.

Observe that one can rewrite the presentation so that the elements
$\varepsilon_{ij}$ with $i,j=0,1$ generate the group.

The case $n=3$ provides another presentation of the fundamental group of a 
sextic with $9$ cusps, different to the one given by Zariski \cite{zar:genusp}
or Kaneko \cite{kan}. In fact they calculate the fundamental group as 
the kernel of the projection of the braid group with three strings on 
the torus, $B_3(T)$, onto the first homology group $H_1(T,{\bf Z})$. The 
presentation had generators 
$$g_2, \ g_{00}, \ g_{01}, \ g_{10} \ \ {\rm and}\ \ g_{11}$$ 
and relations
$$g_2 \cdot g_{ij} \cdot g_2 = g_{ij} \cdot g_2 \cdot g_{ij} \  \ 
{\rm for} \  \ i,j=0,1,2,$$
where 
$$\array{l}
g_{20}=g_{10} \cdot g_{00} \cdot g_{10}^{-1} \\
g_{21}=g_{11} \cdot g_{10} \cdot g_{11}^{-1} \\
g_{02}=g_{01} \cdot g_{00} \cdot g_{01}^{-1} \\
g_{12}=g_{11} \cdot g_{10} \cdot g_{11}^{-1} \\
g_{22}=g_{21} \cdot g_{20} \cdot g_{21}^{-1} \ .
\endarray$$

An isomorphism between both presentations can be given by
$$\array{rll}
\varepsilon_{00} & \mapsto & g_2 \cdot g_{11} \cdot g_2^{-1} \\
\varepsilon_{10} & \mapsto & g_2 \\
\varepsilon_{01} & \mapsto & g_{10} \\
\varepsilon_{11} & \mapsto & g_2 \cdot g_{01} \cdot g_2^{-1} \\
\varepsilon_{11} \cdot \varepsilon_{00}^{01} \cdot \varepsilon_{11}^{-1} 
& \mapsto & g_{00} \ .\\
\endarray$$

\section{Alexander polynomial of $\tilde{C}_n$ for $n$ odd.}
By Theorem $[5.1.]$ in \cite{libg:alex}, we only have to calculate the 
superabundance 
$$s_{\kappa}(2n-3-2n \kappa)=\dim H^1({\cal A}_{\kappa}(2n-3-2n\kappa)),$$
where $\kappa$ is a constant of quasi-adjunction corresponding to a
singularity of type $A_{n-1}$, say $\kappa = \frac{n-2}{2n}$. Observe 
that, for this constant, the local ideal of quasi-adjunction 
$({\cal A}_{\kappa})_{m_p}$ imposes the weakest conditions on the 
linear system, in fact, it is just the maximal ideal of the local ring
at the singular points $p$ of $\tilde{C}_n$.

On one hand, the expected dimension of the linear system 
of curves of degree $n-1$ passing through $3n$ points in general position is 
$$\frac{n(n+1)}{2}-3n.$$
On the other hand, since all singularities belong to a cubic, the dimension 
of the linear system $\Gamma({\bf P}^2,{\cal A}_{\kappa}(n-1))$ is
$$\frac{(n-3)(n-2)}{2}.$$
Therefore,
$$s_{\kappa}(2n-3-2n \kappa )=s_{\frac{n-2}{2n}}(n-1)=
\frac{(n-3)(n-2)}{2}-\Big( \frac{n(n+1)}{2}-3n \Big) = 3.$$
Hence the Alexander polynomial of $\tilde{C}_n$ is

\begin{equation}
\label{alex}
\Delta_{\tilde{C}_n}(t)=(t^{n-1}-t^{n-2}+...+t^2-t+1)^3.
\end{equation}
One can also use (\ref{pres}) to give a presentation of the commutator
and calculate the exponent of (\ref{alex}) as the rank of its 
abelianization.

\section{A deformation of the family}
Instead of considering the conic and the three tangent lines, we will move
one of the lines so that it becomes transverse. There is a surjective map

\begin{equation}
\label{sur1}
\pi_1({\bf P}^2 \setminus (C \cup L_1 \cup L_2 \cup L_3)) \longrightarrow
\pi_1({\bf P}^2 \setminus (C \cup L_1 \cup L_2 \cup L'_3)).
\end{equation}
The same Hirano construction as before will give us a family of curves
$\tilde{C}'_n$ of degree $2n$ but with $2n$ singularities of type $A_{n-1}$
instead of $3n$. We will have a surjection 
\begin{equation}
\label{sur2}
\pi_1({\bf P}^2 \setminus \tilde{C}_n) \longrightarrow
\pi_1({\bf P}^2 \setminus \tilde{C}'_n).
\end{equation}
One can write down the map (\ref{sur1}) and apply the Reidemeister-Schreier 
method simultaneously to get the surjective map (\ref{sur2}). One sees that 
it comes from adding the relations
$$\varepsilon_{i,j}=\varepsilon_i \ {\rm \ for \ any }\  j=0,...,n-1.$$
That is, the non-abelian group 
$({\bf Z}/2{\bf Z}) \ast ({\bf Z}/n{\bf Z})$ 
given in \cite{oka:nonab}. Actually the curves $C_{2,n}$ presented by Oka 
are a deformation of the $\tilde{C}'_n$. If we consider the conic 
$C$ of equation $\{ (x+y+z)^2 + xy =0 \}$ and the abelian covering ramified 
on the coordinate axis, then $\tilde{C}_n$ turns out to be 
$\{ (x^n+y^n+z^n)^2 + x^ny^n =0 \}$, which can be deformed into 
$C_{2,n} \equiv \{ (y^n-z^n)^2 + (x^2-y^2)^n \}$.

{\small Department of Mathematics, Statistics and Computer Science\\ 
University of Illinois at Chicago\\ 851 S. Morgan Street\\ 
Chicago, Illinois, 60607-7045\\ USA\\ e-mail: jicogo@math.uic.edu}

\vspace*{12pt}

{\small Departamento de \'Algebra\\ 
Universidad Complutense de Madrid\\ 
Av. Ciudad Universitaria s/n\\ 28040 Madrid\\ Espa\~na\\
e-mail: jicogo@eucmos.sim.ucm.es}


\begin{thebibliography}{11}

\bibitem{hirano}
A. Hirano,
`Construction of plane curves with cusps',
{\em Saitama Math. J.} 10 (1992) 21--24.
%
\bibitem{kan}
J. Kaneko,
`On the fundamental group of the complement to a maximal 
cuspidal plane curve',
{\em Memoirs of the Fac. of Sc. Kyushu University},
Ser.~A, 39, No.~1 (1985) 133--146.
%
\bibitem{libg:alex}
A. Libgober,
`Alexander invariants of plane algebraic curves',
{\em Proc. of Symp. in Pure Math.}
40 (1983) Part~2, 135--143.
%
\bibitem{libg:homot}
A. Libgober, 
`Homotopy groups of the complements to singular 
hypersurfaces. II', 
{\em Ann. of Math.} (2)
139 (1994) No.~1, 117--144.
%
\bibitem{kar:comb}
W. Magnus,--A. Karras,--D. Solitar,
{\em Combinatorial group theory}
(New York, Interscience Publishers, 1966).
%
\bibitem{mois:braid}
B. Moishezon,
`Stable branch curves and braid monodromies',
{\em Lectures Notes in Math.}
862, Springer-Verlag, Berlin and New York, 1981.
%
\bibitem{oka:nonab}
M. Oka,
`Some plane curves whose complement have non-abelian fundamental groups',
{\em Mat. Ann.}
218 (1978) 55--65.
%
\bibitem{oka:two}
M. Oka,
`Two transforms of plane curves and their fundamental groups',
{\em J. Math. Sci. Univ. Tokio}
3 (1996) 399--443.
%
\bibitem{zaid}
S. Orevkov, and M. Zaidenberg,
`On the number of singular points of plane curves',
{\em Duke preprint alg-geom/9507005}.
%
\bibitem{shim}
I. Shimada,
`A note on Zariski pairs',
{\em Comp. Math.}
104 (1996) No.~2, 125--133.
%
\bibitem{zar:genusp}
O. Zariski,
`The topological discriminant group of a Riemann surface of genus p',
{\em Am. J. of Math.}
59 (1937) 335--358.
%
\bibitem{zar2}
O. Zariski,
`On the problem of existence of algebraic functions of two variables
possesing a given branch curve',
{\em Am. J. of Math.},
51 (1929) 305--328.

\end{thebibliography}
\end{document}